\documentclass[11pt]{amsart}
\usepackage{amsfonts,amsmath,amssymb}
\usepackage[all]{xy}
\usepackage{hyperref}

\begin{document}

\parindent=0pt
\parskip=6pt


\newcommand{\F}{{\mathbb F}}
\newcommand{\C}{{\mathbb C}}
\newcommand{\Q}{{\mathbb Q}}
\newcommand{\mP}{{\mathbb P}}
\newcommand{\R}{{\mathbb R}}
\newcommand{\T}{{\mathbb T}}
\newcommand{\Z}{{\mathbb Z}}
\newcommand{\Sl}{{{\rm Sl}_2}}
\newcommand{\vk}{{\varkappa}}
\newcommand{\dR}{{\rm dR}}
\newcommand{\be}{{\bf e}}
\newcommand{\bu}{{\bar{u}}}
\newcommand{\bh}{{\mathfrak{h}}}
\newcommand{\bv}{{\bar{v}}}
\newcommand{\bw}{{\bf v}}
\newcommand{\bM}{{\overline{M}}}
\newcommand{\tM}{{\widetilde{M}}}
\newcommand{\PGl}{{{\rm PGl}_2}}
\newcommand{\SU}{{\rm SU}}
\newcommand{\vt}{{\vartheta}}
\newcommand{\tT}{{\tilde{t}}}
\newcommand{\ts}{{\tilde{\sigma}}}
\newcommand{\half}{{\textstyle{\frac{1}{2}}}}
\newcommand{\quarter}{{\textstyle{\frac{1}{4}}}}
\newcommand{\third}{{\textstyle{\frac{1}{3}}}}
\newcommand{\conf}{{\rm Config}}
\newcommand{\sT}{{\sf T}}

\newcommand{\ie}{\textit{i}.\textit{e}\;}
\newcommand{\eg}{\textit{e}.\textit{g}\;}
\newcommand{\cf}{{\textit{cf}\;}}

\title{On oriented planar trees with three leaves}

\author[Jack Morava]{Jack Morava}

\address{Department of Mathematics, The Johns Hopkins University,
Baltimore, Maryland 21218}

\email{jack@math.jhu.edu}


\date{May 2020}

\begin{abstract}{This elementary note proposes candidates for 
interesting continuous piecewise-smooth `Riemannian' metrics 
on the moduli spaces of rooted geodesic trees embedded in the
Poincar\'e disk. A related digression observes the existence 
of an apparently hitherto unrecognized abelian topological 
group structure on the real projective line.

\noindent This is an early draft. The symbol $\star$ indicates
an illustration to be supplied in a later version.}\end{abstract}
\bigskip

\maketitle \bigskip

{\bf 0.1 Introduction} A finite set of points on the boundary of the 
Poinc\'are unit disk can be shown to define a unique embedded geodesic
not necessarily binary tree, with one point chosen as root and the rest 
as endpoints of its leaves. Alternatively, the dual \cite{2}(\S 1) of 
such a spanning tree defines a decomposition of the disk into hyperbolic 
polygons, which can be imagined to be solutions of a geometric extremization 
problem analogous to Plateau's\begin{footnote}{\ie to find a minimal surface 
spanning a given closed space curve.}\end{footnote}, \ie to find the pattern 
of cracks in a shattered windshield, given a collection of shocks or stresses 
at a finite number of boundary points. As a problem in analytic mechanics 
this seems to be a relatively intractable kind of free boundary problem, 
making a useful account of its literature problematic. The present sketch 
is concerned with the differential geometry of the space $\bM_{0,n+1}(\R)$ 
of such rooted geodesic hyperbolic trees, as providing an alternate approach 
to questions of this sort.\bigskip 

{\bf 0.2 $\star$} The simplest nontrivial case (rooted trees with three leaves)
is classical: the cross-ratio $\rho = [x_0:x_1:x_2:x_3]$ maps the space
$\bM_{0,4}(\R)$ (of projective equivalence classes of four labelled points
on the real projective line) diffeomorphically to the circle. When $\rho
\in [0,1]$, a beautiful formula 
\[
\rho = \frac{1}{1 + e^{-\gamma}}
\]
of Devadoss equates the cross-ratio to the logistic function of the signed
hyperbolic length of the generically unique internal branch. 

In general $\bM_{0,n+1}(\R)$ is tesselated by $\half n!$ Stasheff associahedra,
which can be regarded as moduli spaces for geodesic trees with leaves in 
fixed order. This formula then generalizes \cite{3}(\S 6) to define a 
pseudometric which blows up along the associahedral faces. On the other hand
the original construction of $\bM_{0,n+1}(\R)$ (as the real points of the
Deligne-Mumford-Knudsen compactification of the stack of genus zero algebraic
curves with marked points) makes it a smooth projective variety, which thus
inherits a very nice Riemannian metric from its K\"ahler structure - which 
is however not unique, or necessarily well-behaved with respect to permutations
of the marked points. 

{\bf 0.3} The purpose of this note is to suggest that pulling back (a 
variant of) the product metric along a generalized Albanese map
\[
\bM_{0,n+1}(\R) \to \bM_{0,3+1}(\R)^{\binom {n}{3}} 
\]
provides the moduli space with a continuous, though only piecewise smooth, 
$\Sigma_n$-equivariant metric; and, more generally, to suggest the relevance
of the metric geometry of these spaces to the study of such generalized 
travelling salesman problems. \bigskip

{\bf \S I The real projective line, reconsidered}\bigskip

{\bf 1.1} The projective line $\mP(\F_2) = \{0,1,\infty\}$ over the field with 
two elements admits a transitive action of the symmetric group $\Sigma_3 \cong
\Sl(\F_2)$. The rational functions
\[
\tau_{01}(x) = 1 - x, \; \tau_{1 \infty}(x) = (1 - x^{-1})^{-1}
\]
generate a group (under composition) isomorphic to $\Sigma_3$, satisfying the
braid relation 
\[
(\tau_{01} \circ \tau_{1 \infty} \circ \tau_{01})(x) =  (\tau_{1 \infty} \circ
\tau_{01} \circ \tau_{1 \infty})(x) = x^{-1} = 
\left[\begin{array}{cc}
                   0 & 1 \\                        
                   1 & 0 \end{array}\right](x) \;.
\]
The element
\[
\sigma(x) = (\tau_{1\infty} \circ \tau_{01})(x) = (1 - x)^{-1} =
\left[\begin{array}{cc} 
                   0 & 1 \\
                  -1 & 1 \end{array}\right](x)
\]
(sending $(0 1 \infty)$ to $(\infty 0 1)$) generates a cyclic subgroup 
$C_3 = \{1,\sigma,\sigma^2\}$ of $\Sigma_3$, with 
\[
\tau_{01} \circ \tau_{1\infty} = \sigma^2(x) = 1 - x^{-1} = 
\left[\begin{array}{cc}
                  1 & -1 \\
                  1 &  0 \end{array}\right](x) \;;
\]
thus $\tau_{01} \circ \sigma \circ \tau_{01} = \sigma^{-1}$. This lifts to an 
action of $C_3 \rtimes C_2 \cong \Sigma_3$ on $\C$ by rational functions, 
which extends to a continuous (fixed-point free) action on the one-point 
compactification $\R_+ = \R \cup \infty = \mP_1(\R)$ of the real line. 
However $\sigma(x)$ is not differentiable at $x = 0$, and $\sigma'(x) = 0$ 
at $x = \infty$. \bigskip

{\bf 1.2} Let $(x_i) \in \R^4, \; 0 \leq i \leq 3$, be a vector with all 
coordinates distinct, and define $x_{ij} = x_i - x_j \neq 0$ for $i \neq j$; 
then the cross-ratio
\[
[x_0:x_1:x_2:x_3] = \frac{x_{01}}{x_{02}} \cdot \frac{x_{23}}{x_{13}}
\]
extends to the compactified quotient
\[
\xymatrix{
{\conf^4(\R)/\PGl(\R) \subset \bM_{0,4}}(\R) \ar[r]^-\rho_-\cong
& \mP_1(\R) } \;,
\]
defining an isomorphism of the space of four points on the line, modulo
projective equivalence, with the projective line itself. Thus the 
fractional linear transformation 
\[
x \mapsto [A](x) = \frac{ax + b}{cx + d}
\]
defined by 
\[
A = \left[\begin{array}{cc}
                   a & b \\
                   c & d \end{array}\right] \in \PGl(\R)
\]
(\ie $ad - bc = \det T \neq 0$) satisfies
\[
[[A](x_0):[A](x_1):[A](x_2):[A](x_3)] = [x_0:x_1:x_2:x_3] \;;
\]
for example
\[
\rho = [0:\rho:1:\infty] = [1:\sigma(\rho):\infty:0] = 
[\infty:\sigma^2(\rho):0:1] \;,
\]
\eg $\rho = \sigma(\rho)^{-1}(\sigma(\rho) - 1)$. 

{\bf Remark} If $\mP_1(\R)$ is ordered as usual (\ie with $x \in \R \Rightarrow
x < \infty$), then for any quadruple $x_0 < x_1 < x_2 < x_3$ there is an
$A \in \Sl(\R)$ such that $[A](x)_0) = 0, \; [A](x_2) = 1, \; [A](x_3) =  
\infty$, and $[A](x_1) \in (0,1)$. If $x_3 \neq \infty$, then
\[
A = (x_{01}x_{20}x_{23})^{-1/2}
\left[\begin{array}{cc}
             x_{23} & -x_1x_{23} \\
             x_{20} & -x_0x_{20} \end{array}\right]
\]
is unique.

The group
\[
\Sigma_4 \cong \F^2_2 \rtimes \Sl(\F_2)
\]
of permutations of four things acts on $\bM_{0,4}(\R)$, but the Klein 
subgroup $\F^2_2$ generated by the transpositions
\[
(01)(23), \; (02)(13), \; (12)(30)
\]
leaves the cross-ratio invariant, so the action of $\Sigma_4$ on 
$\bM_{0,4}(\R)$ reduces to the action of $\Sigma_3$ described above.
\bigskip

{\bf 1.3} $\star$  It will be convenient to decompose the projective line
\[
\mP_1(\R) = [I] \cup [II] \cup [III]
\]
as the union of closed intervals 
\[
[I] = [-\infty,0],\; [II] = [0,1],\; [III] = [1,\infty] \;;
\]
thus $\sigma$ maps $[I]$ bijectively to $[II]$, $[II]$ to $[III]$, and
$[III]$ to $[I]$, preserving orientations, and is consistent on the 
boundary points 
\[
[I] \cap [II] = \{0\}, \; [II] \cap [III] = \{1\}, \; [III] \cap [I] = 
\{\infty\} \;.
\]
Similarly, let $(X)$ denote the interior of $[X]$; then the (smooth) function 
\[
k : \R - \{0,1\} = (I) \cup (II) \cup (III) \to (0,1)
\]
defined by $k(x) =$
\[
\sigma(x) = (1 - x)^{-1} \; {\rm if} \; x \in (-\infty,0) = (I) \;,
\]
\[
= x \; {\rm} \; {\rm if} \; x \in (0,1) = (II) \;, \; {\rm and}
\]
\[
= \sigma^{-1}(x) = 1 - x^{-1} \; {\rm if} \; x \in (-\infty,0) = (III) 
\]
extends to a continuous function 
\[
\kappa : \R_+ = \mP_1(\R) \to (0,1)_+ = \R/\Z
\]
sending $\{0,1,\infty\}$ to the compactification point $0 = 1 \in (0,1)_+$.  
\bigskip

{\bf Lemma} $\star$ {\it The derivative $\kappa'(x)$
\[
= (1 - x)^{-2}, \; x \in (I)
\]
\[
= 1, \; x \in (II)
\]
\[
= x^{-2}, \; x \in (III) 
\]
exists and is continuous. Moreover,} 
\[
\int_\R \kappa'(x) \cdot dx = \int_{-\infty}^0 (1 - x)^{-2} \cdot dx + 
\int_0^1 1 \cdot dx + \int_1^\infty x^{-2} \cdot dx  = 1+1+1 = 3 \;.
\]

{\bf Proposition} $\kappa \circ \sigma = \kappa$: {\it thus $\kappa$ identifies the
quotient of $\mP_1(\R)$ by $C_3$ with $\R/\Z$.}

[For if $x \in (I)$ then $\sigma(x) \in (II)$ so $\kappa (\sigma(x)) = 
\sigma(x) = \kappa(x)$; while if $x \in (II)$ then $\sigma(x) \in (III)$,
so $\kappa(\sigma(x)) = \sigma^{-1}(\sigma(x)) = x = \kappa(x)$. Finally, if 
$x \in (III)$ then $\sigma(x) \in (I)$ so 
\[
\kappa(\sigma(x)) = (1 - \sigma(x))^{-1} = 1 - x^{-1} = \kappa(x) \;.]
\]

The resulting map is a three-fold cover, with multiplication by 3 as the
induced homomorphism
\[
\pi_1(\mP_1(\R),\infty) \cong \Z \to \Z \cong \pi_1(\R/\Z,0) \;.
\]
The action of $\Sigma_3$ on $\mP_1(\R)$ reduces to the orientation-reversing
action of $\Sigma_3/C_3 \cong \{\pm 1\}$ on $\mP_1(\R)/C_3$ defined by 
$\tau_{0,1}(x) = 1 - x$ on $\R/\Z$ as in \S 1.2. The composition
\[
\xymatrix{
\bM_{0,3+1}(\R) \ar[r]^-\rho & \mP_1(\R) \ar[r]^-\kappa & \R/\Z }
\]
provides an interpretation the quotient space $\bM_{0,3+1}(\R)/C_3$
as the space of configurations $(\{x_i,x_j,x_k\},\infty)$ of three 
cyclically ordered points on $\R$. The one-form $d\left\lceil t 
\right\rceil \in \Omega^1(\R/\Z)$ pulls back to a $C_3$-invariant one-form
\[
d\kappa = \kappa'(\rho)d\rho = \kappa^*(d\rho) \in \Omega^1(\bM_{0,3+1}(\R))
\]
mapping to three times the fundamental class in $H^1_\dR(\bM_{0,3+1}(\R))$.
\bigskip

{\bf 1.4} $\star$ The inverse 
\[
D(x) = - \log |1 - x^{-1}| = - \log |\sigma^{-1}(x) : [0,1] \to \mP_1(\R)
\]
of the logistic function $(1 + e^{-x})^{-1}$ extends to the closed interval
by $D(0) = D(1) = \infty$. Thus $\vk (x) = (D \circ \kappa)(x)$ 
\[
= - \log |x| \; {\rm if} \; x \in [I] \;, 
\]
\[
= - \log |\sigma^{-1}(x)| \; {\rm if} \; x \in [II] \;, 
\]
\[ 
= - \log |\sigma^{-2}(x)| \; {\rm if} \; x \in [III]  
\]
defines a three-fold cover $\vk : (\mP_1(\R),\{0,1,\infty\}) \to (\mP_1(\R),
\infty)$ of the projective line. Its  graph
\[
\mP_1(\R) \ni x\mapsto (x,\vk(x)) \in \mP_1(\R) \times \mP_1(\R)
\]
has degree one along the first factor, and degree three
along the second, defining a  piecewise smooth helix wrapped 
around a torus. 

{\bf Proposition} $\star$  
\[
\gamma = \vk(\rho) : \bM_{0,3+1}(\R) \to \mP_1(\R)
\]
{\it extends Devadoss's formula \cite{2}(\S 6.1) for the oriented hyperbolic 
length of the generic internal edge of a hyperbolic rooted tree with three 
leaves. The square of the one-form 
\[
d\gamma = \vk'(\rho)d\rho
\]
defines a pseudometric on $\bM_{0,3+1}(\R)$ which blows up at} $\rho = 0,1,
\infty$.\bigskip

{\bf 1.5} The rational cohomology $H^*(\bM_{0,n+1}(\R),\Q)$ of the moduli
space of $n+1$ ordered points on the line is calculated in \cite{4}(\S 2.3),
\cite{7}; in particular
\[
H^1(\bM_{0,n+1}(\R),\Q) \cong \Lambda^3 \bh_n
\]
as $\Sigma_{n+1}$-modules, where $\bh_n$ is the $n$-dimensional kernel of the 
trace homomorphism
\[
\Q^{n+1} \ni (v_0,\dots,v_n) \to \sum v_i \in \Q
\]
(with $\Sigma_{n+1}$ acting on the left by permuting coordinates). A subset
$S$ of $\{1,\dots,n\}$ of cardinality $3 \leq |S| \leq n$ defines a forgetful
morphism $\bM_{0,n+1}(\R) \to \bM_{0,|S|+1}(\R)$: thus a subset
$S = \{i < j < k\}$ defines a composition
\[
\xymatrix{
\bM_{0,n+1}(\R) \ar[r]^-{\nu_S} & \bM_{0,3+1}(\R) \ar[r]^-\kappa & \R/\Z } \;.
\]
and hence  a one-form $d\kappa_S \in \Omega^1(\bM_{0,n+1}(\R))$. I will write
$d\kappa_S^{\otimes 2} \in \Omega^{\otimes 2}(\bM_{0,n+1}(\R))$ for the 
associated quadratic differential. 

A basis $\bw_i$ for $\Q^{n+1}$ defines a basis 
\[
\alpha_{ijk} = (\bw_i - \bw_0) \wedge (\bw_j - \bw_0) \wedge (\bw_k - \bw_0)
\]
for $\Lambda^3\bh_n$; then 
\[
\Lambda^3 \bh_n \ni \alpha_{ijk} \mapsto d\kappa_{ijk} \in \Omega^1
(\bM_{0,n+1}(\R))
\]
is an injective homomorphism of $\Sigma_n$-modules. 

{\bf Claim} The average
\[
ds^2_\bM := \binom{n}{3}^{-1} \sum_{i<j<k} d\kappa_{ijk}^{\otimes 2} \in 
\Omega^{\otimes 2} (\bM_{0,n+1}(\R))
\]
{\it defines  a continuous, piecewise smooth metric on the space of rooted 
hyperbolic $n$-leaved trees. The subgroup $\Sigma_n \subset \Sigma_{n+1}$ of 
permutations acts by isometries.}

[Behind this lies the conjecture that the `Albanese' map
\[
\prod_{S \in \binom{n}{3}} \kappa_S : \bM_{0,n+1}(\R) \to 
\bM_{0,3+1}(\R)^{\binom{n}{3}}
\]
is an immersion.]\bigskip 

{\bf \S II The Cayley transform, reconsidered}\bigskip

{\bf 2.1} The fractional linear transformation
\[
z = C(x) = 
\left[\begin{array}{cc}
                   1 & -i \\
                  - i & 1 \end{array}\right](x) =
\frac{x - i}{1 - ix} : \mP_1(\C) \to \mP_1(\C)
\]
restricts to stereographic projection 
\[
\mP_1(\R) \supset \R \to \T \subset \C^\times \subset \mP_1(\C)
\]
of the real line, sending $\pm \infty \to i$ and $\pm 1$ to $\pm 1$. Writing
$t \mapsto \be(t) = \exp(2 \pi it)$ for the group isomorphism 
$\R/\Z \to \T$, we have
\[
\left[\begin{array}{cc}
                   1 & -1 \\
                   1 & 1 \end{array}\right](\be(t)) = i \tan \pi t \;,
\]
so (by the addition formula for the tangent function)
\[
x = C^{-1}(\be(t)) =
\left[\begin{array}{cc}
                   1 & i \\
                   i & 1 \end{array}\right] \cdot
\left[\begin{array}{cc}
                   1 & 1 \\
                 - 1 & 1 \end{array}\right](i \tan \pi t) =
\left[\begin{array}{cc}
                   1-i & 1+i \\
                   i-1 & 1-i \end{array}\right](i \tan \pi t) 
\]
\[
= \frac{1 + \tan \pi t}{1 - \tan \pi t} = \tan \pi (t + \quarter)
\]
defines a diffeomorphism (inverse to stereographic projection) of $\T$ with 
$\R_+ = \mP_1(\R)$. Similarly, $(2\pi iz)^{-1}dz \in \Omega^1(\C^\times)$
pulls back to $(\pi (1 + x^2))^{-1}dx \in \Omega^1(\R)$. \bigskip

{\bf 2.2} The action of $A \in \Sl(\R)$ on $\mP_1(\R)$ defines an action of 
\[
\tilde{A} = C^{-1}AC = (\half)^2
\left[\begin{array}{cc}
                   1 & i \\
                   i & 1 \end{array}\right] \cdot
\left[\begin{array}{cc}
                   a & b \\
                   c & d \end{array}\right] \cdot
\left[\begin{array}{cc}
                   1 & -i \\
                  -i & 1 \end{array}\right]
\]
\[
= \left[\begin{array}{cc}
                   u & v \\
                 \bv & \bu \end{array}\right] \in \SU(1,1) \;,
\]
where $|u|^2 - |v|^2 = 1$, with
\[
u = \half [(a+d) + i(c-b)], \; v = \half[(b+c) + i(d-a)] \;,
\]
on $\C$. Because the complex conjugate of 
\[
\left[\begin{array}{cc}
                   u & v \\
                 \bv & \bu \end{array}\right](\be(t))
\]
equals its inverse
\[
\frac{\bu \be(-t) + \bv}{v \be(-t) + u} = \frac{\bv \be(t) + \bu}
{u \be(t) + v} \;,
\]
this action takes the circle $\T \subset \C^\times$ to itself. 

{\bf Example} $\star$  
\[
\ts = \half \left[\begin{array}{cc}
                   1-2i & +i \\
                     -i & 1+2i \end{array}\right], \;
\ts^2 = \half \left[\begin{array}{cc}
                   1+2i & -i \\
                     +i & 1-2i \end{array}\right],
\]
so $\{1,\pm i\}$ is an orbit of $\tilde{C}_3$, and $\{-1,(4 \pm 3i)/5\}$
is another. Similarly, 
\[
\tilde{\tau}_{01} = 
\frac{i}{2} \left[\begin{array}{cc}
                   -i & 1 + 2i \\
                   1 - 2i & i \end{array}\right] \;.
\]\bigskip

{\bf 2.3} The (renormalized) extension 
\[
z \mapsto L(z) = iC(z) = \frac{1 + iz}{1 - iz} =
\left[\begin{array}{cc}
                   i & 1\\
                  -i & 1\end{array}\right](z)
\]
of $C$ to a fractional linear transformation $L$ of the complex projective 
line $\mP_1(\C) \cong \C_+$ maps $\mP_1(\R)$ to $\T$. .

{\bf Proposition}
\[
L^{-1}(L(z_0) \cdot L(z_1))\; = \; \frac{z_0 + z_1}{1 - z_0z_1} := z_0 +_L z_1
\]
{\it restricts near 0 to the one-dimensional formal group law with $ x \mapsto
\tan x$ as its exponential.} 

{\bf Proof}
\[
L^{-1}(z) =
\left[\begin{array}{cc}
                   1 & -1\\
                   i & i\end{array}\right](z) = -i \frac{z-1}{z+1} \;,
\]
so
\[L^{-1}(L(z_0) \cdot L(z_1)) =
L^{-1}[\frac{1 + iz_0}{1 - iz_0} \cdot \frac{1 + iz_1}{1 - iz_1}] =
\]
\[
(-i) \cdot \frac{(1 + iz_0)(1 + iz_1) - (1 - iz_0)(1 - iz_1)}{(1 + iz_0)
(1 + iz_1) + (1 - iz_0)(1 - iz_1)} =
\]
\[
(-i) \frac{2i(z_0 + z_1)}{2 - 2z_0z_1} = z_0 +_L z_1
\]
as claimed. $\Box$

Up to a Wick rotation $x \mapsto ix$, this is the formal group law defined by
Weyl and Hirzebruch's signature genus for oriented smooth manifolds.  It is
odd, in the sense that $[-1]_L(z) = -z$.

{\bf Corollary} $(\mP_1(\R) = \R_+,0,+_L)$ {\it  is a group, with $\infty$ 
as (the unique nontrivial) torsion point of order two}.

In particular, if $x \in \R^\times$, then
\[
x \mapsto x +_L \infty = \lim_{w \to 0} \frac{x + w^{-1}}{1 - w^{-1}x} =
\lim_{w \to 0} \frac{1 + wx}{w - x} = -x^{-1} \;.
\]
Similarly, $x +_L x = [2]_L(x) \to 0$ as $x \to \infty$, consistent with
$\infty = [-1]_L(\infty)$, while
\[
[3](w^{-1}) = \frac{3w^2 - 1}{w^3 - 3w} \to \infty
\]
as $w \to 0$, etc.

Note also that 1 is a 4-torsion point, \ie $1 +_L 1 = \infty$. More generally,
the group $\Q/\Z \subset \T$ of torsion points for $+_L$ maps isomorphically
to the set $\{ \tan \pi x \: | \: x \in \Q \}$ of (cyclotomic) algebraic 
numbers.  

{\bf Note} that since $i +_L (-i)$ is undefined, this construction fails to
define a composition operation on  $\mP_1(\C)$. \bigskip

\bibliographystyle{amsplain}

\end{document}